\documentclass[11pt]{article}

\usepackage{amssymb}

\newcommand{\tends}[1]{{\displaystyle\mathop{\to}_{#1}}}

\def\C{{\mathbb C}}
\def\N{{\mathbb N}}

\def\R{{\mathbb R}}
\def\T{{\mathbb T}}
\def\Z{{\mathbb Z}}

\def\Ker{{\rm Ker}\,}
\def\Mat{{\rm Mat}}
\def\SL{{\rm SL}}
\def\Spec{{\rm Spec}\,}

\def\eps{\varepsilon}

\title{\bf On the K-property of quantized Arnold cat maps}

\author{S.V. Neshveyev}

\date{}

\begin{document}

\maketitle

\begin{abstract}
We prove that some quantized Arnold cat maps are entropic K-systems. This
result was formulated by H.~Narnhofer \cite{Nar}, but the fact that the
optimal decomposition for the multi-channel entropy constructed there is
not strictly local was not appropriately taken care of. We propose a
strictly local decomposition based on a construction of Voiculescu.
\end{abstract}

\bigskip

\section{Introduction}
The concept of K-system is very important in ergodic theory. Narnhofer and
Thirring \cite{NT} introduced a non-commutative analogue
of this notion. In \cite{GN1} V.Ya.~Golodets and the author proved the
following sufficient condition for the K-property: a W$^*$-system
$(M,\phi,\alpha)$ is an entropic K-system if there exists a W$^*$-subalgebra
$M_0$ of $M$ such that $M_0\subset\alpha(M_0)$,
$\cap_{n\in\Z}\alpha^n(M_0)=\C1$,
$\cup_{n\in\N}(\alpha^{-n}(M_0)'\cap\alpha^n(M_0))$ is weakly dense in $M$.
This condition and the observation that a subsystem of a K-system invariant
under the modular group is a K-system too allow to construct a large
class of quantum K-systems (see, in particular, \cite{GN1,GN2,N}). We know
only one class of
quantum systems for which the K-property is obtained by different
arguments. This is quantized Arnold cat maps. This result was formulated in
Narnhofer's paper \cite{Nar}. The decompositions contructed in the
course of the proof there are not strictly local, that leads to a
factor that again could only be controlled by using asymptotic
abelian arguments. So the essential interest lies in the
construction of a completely positive map that is strictly local
and can be well controlled and generalized in a larger context.

\bigskip

\section{The K-property of quantized cat maps}
Let $G$ be a discrete abelian group, $\omega\colon G\times G\to\T$ a
bicharacter. Consider the twisted group C$^*$-algebra $C^*(G,\omega)$
generated by unitaries $u_g$, $g\in G$, such that
$$
u_gu_h=\omega(g,h)u_{g+h}.
$$
The canonical trace $\tau$ on $C^*(G,\omega)$ is given by $\tau(u_g)=0$ for
$g\ne0$. It is known that the uniqueness of the trace is equivalent to the
simplicity of $C^*(G,\omega)$, and is also equivalent to the non-degeneracy
of the pairing $(g,h)\mapsto\omega(g,h)\bar\omega(h,g)$. In particular, if
$G$ is countable and the pairing is non-degenerate, then
$\pi_\tau(C^*(G,\omega))''$ is the hyperfinite II$_1$-factor. Each
$\omega$-preserving automorphism $T$ of $G$ defines an automorphism
$\alpha_T$ of $C^*(G,\omega)$, $\alpha_T(u_g)=u_{Tg}$.

The non-commutative torus $A_\theta$ ($\theta\in[0,1)$) is the algebra
$C^*(\Z^2,\omega_\theta)$, where
$$
 \omega_\theta(g,h)=e^{i\pi\theta\sigma(g,h)},\ \sigma(g,h)=g_1h_2-g_2h_1.
$$

The following theorem was formulated in \cite{Nar}.

\medskip\noindent
{\bf Theorem 1.} {\it Let $T\in\SL_2(\Z)$, $\Spec T=\{\lambda,\lambda^{-1}\}$.
 Suppose $|\lambda|>1$ (so that $\lambda$ is real) and
 $\theta\in[0,1)\cap(2\Z\lambda^2+2\Z)$. Then
 $(\pi_\tau(A_\theta)'',\tau,\alpha_T)$ is an entropic K-system.}

\medskip
We will prove the following more general result.

\medskip\noindent
{\bf Theorem 2.} {\it Let $T$ be an aperiodic $\omega$-preserving
 automorphism of $G$. Suppose that
 $$
 \sum_{n\in\Z}|1-\omega(g,T^nh)|<\infty\ \ \forall g,h\in G.
 $$
 Then $(\pi_\tau(C^*(G,\omega))'',\tau,\alpha_T)$ is an entropic K-system.}

\medskip
It was proved in \cite[Theorem 3.8]{Nar} that under the assumptions of
Theorem 1, for any $g,h\in\Z^2$, we have
$$
|1-\omega(g,T^nh)|\le C|\lambda|^{-|n|},
$$
so Theorem 1 is really follows from Theorem 2. The key observation for
that estimate was the equality
$$
T^nh={1\over\lambda^2-1}\sum^2_{i=0}(\lambda^{n+i}+\lambda^{-n-i})\bar h_i
       +\lambda^{-n}\bar h\ \ (n\in\N),
$$
where $\bar h_i\in\Z^2$ and $\bar h\in\R^2$ depend only on $h$ and $T$,
which is obtained by computations in a basis diagonalizing $T$. Since
$\sigma(g,\bar h_i)$, $\sigma(g,h)$,
$\lambda^{n+i}+\lambda^{-n-i}=\hbox{Tr}\,T^{n+i}$ are all integers and
$\theta\equiv2s(\lambda^2-1)\,\hbox{mod}\,2\Z$ for some $s\in\Z$, we have
$$
\theta\sigma(g,T^nh)\equiv\lambda^{-n}2s(\lambda^2-1)\sigma(g,\bar h)
 \,\hbox{mod}\,2\Z,
$$
whence
$|1-\omega(g,T^nh)|\le|\lambda|^{-n}2\pi|s|(\lambda^2-1)|\sigma(g,\bar h)|$.

\medskip
Starting the proof of Theorem 2, consider a unital completely positive
mapping $\gamma\colon A\to\pi_\tau(C^*(G,\omega)''$ of
a finite-dimensional C$^*$-algebra $A$. By definition \cite{NT}, we have to
prove that
$$
\lim_{n\to\infty}\lim_{k\to\infty}
 {1\over k}H_\tau(\gamma,\alpha^n_T\circ\gamma,
                    \ldots,\alpha^{n(k-1)}_T\circ\gamma)=H_\tau(\gamma).
$$

For a finite set $X$, we denote by $\Mat(X)$ the C$^*$-algebra of linear
operators on $l^2(X)$. Let $\{e_{xy}\}_{x,y\in X}$ be the canonical system
of matrix units in $\Mat(X)$. For $X\subset G$, we define a unital
completely positive mapping $i_X\colon\Mat(X)\to C^*(G,\omega)$ by
$$
i_X(e_{xy})={1\over|X|}u_xu^*_y={\bar\omega(x-y,y)\over|X|}u_{x-y}.
$$
As follows from \cite{V} (see Lemmas 5.1 and 6.1 there), there exist a net
$\{X_i\}_i$ of finite subsets in
$G$ and, for each $i$, a unital completely positive mapping
$j_{X_i}\colon C^*(G,\omega)\to\Mat(X_i)$ such that
$||(i_{X_i}\circ j_{X_i})(a)-a||\tends{i}0$ $\forall a\in C^*(G,\omega)$.
>From this we may conclude that any partition of unit in
$\pi_\tau(C^*(G,\omega))''$ can be approximated in strong operator topology
by a partition of the form $\{i_X(a_k)\}_k$, where $\{a_k\}_k$ is a
partition of unit in $\Mat(X)$. Hence, for any $\eps>0$, there exist
a finite subset $X\subset G$ and a finite partition of unit
$1=\sum_{i\in I}a_i$ in $\Mat(X)$ such that, for $b_i=i_X(a_i)$, we have
$$
H_\tau(\gamma)<\eps+\sum_i\eta\tau(b_i)+
   \sum_iS(\tau(\gamma(\cdot)),\tau(\gamma(\cdot)b_i)),
$$
where $\eta x=-x\log x$.
Set $X_{nk}=\sum^k_{l=1}T^{n(l-1)}(X)$.

The following lemma was proved in \cite{Nar} for $G=\Z^2$.

\medskip\noindent
{\bf Lemma.} {\it Let $G$ be a discrete abelian group, $T$ an aperiodic
 endomorphism of $G$, $\Ker T=0$, $Y$ a finite subset of $G$, $0\in Y$.
 Then there exists $n_0\in\N$ such that if
 \begin{equation} \label{e1}
 \sum^k_{l=1}T^{n(l-1)}y_l=0
 \end{equation}
 for some $y_1,\ldots,y_k\in Y$, $n\ge n_0$, $k\in\N$, then
 $y_1=\ldots=y_k=0$.}

\medskip\noindent{\it Proof.}
First consider the case where $G$ is finitely generated. Then the periodic
part of $G$ is finite. Since $T$ acts on it aperiodically, it is trivial,
so $G\cong\Z^n$ for some $n\in\N$. Then $T$ is defined by a non-degenerate
matrix with integral entries, which we denote by the same letter $T$. It is
known that the aperiodicity is equivalent to $\T\cap\Spec T=\emptyset$.
Let $\Spec T=\{\lambda_1,\ldots,\lambda_m\}$, $V_i\subset\C^n$ be the root
space corresponding to $\lambda_i$, and $P_i$ the projection onto
$V_i$ along $\oplus_{j\ne i}V_j$. Then (\ref{e1}) is equivalent to the
system of equalities
\begin{equation} \label{e2}
\sum^k_{l=1}T^{n(l-1)}P_iy_l=0,
\end{equation}
$i=1,\ldots,m$. Fix $i$. Suppose, for definiteness, that $|\lambda_i|<1$, and
choose $\delta$, $0<\delta<1-|\lambda_i|$. Since $T|_{V_i}$ is a sum of
Jordan cells, there exists a constant $C$ such that
$$
||T^n|_{V_i}||\le C(|\lambda_i|+\delta)^n \ \ \forall n\in\N.
$$
There exists also a constant $M>0$ such that, for $y\in Y$, we have either
$P_iy=0$ or $M^{-1}\le||P_iy||\le M$. Finally, choose $n_i\in\N$ such that
$$
\sum^\infty_{n=n_i}MC(|\lambda_i|+\delta)^n<M^{-1}.
$$
Then if the equality (\ref{e2}) holds with $n\ge n_i$, then $P_iy_1=0$.
Since $\Ker T=0$, we can rewrite (\ref{e2}) as
$\sum^{k-1}_{l=1}T^{n(l-1)}P_iy_{l+1}=0$. Thus we sequentially obtain
$P_iy_1=\ldots=P_iy_k=0$. So we may take $n_0=\max_in_i$.

We prove the general case by induction on $|Y|$ using the same method as
in \cite{R} to reduce the proof to the case considered above.

Let $H_0$ be the group generated by $Y,TY,T^2Y,\ldots\ $. Set $H_n=T^nH_0$,
$H_\infty=\cap_nH_n$, $Y'=Y\cap H_\infty$. Suppose $Y'\ne Y$. There exists
$n_1\in\N$ such that $Y'=Y\cap H_{n_1}$. If the equality (\ref{e1}) holds
with $n\ge n_1$, then $y_1\in H_{n_1}\cap Y=Y'\subset H_\infty$. Then
$\sum^k_{l=2}T^{n(l-1)}y_l\in H_\infty$. Since $\Ker T=0$ and
$TH_\infty=H_\infty$, we conclude that
$\sum^{k-1}_{l=1}T^{n(l-1)}y_{l+1}\in H_\infty$. Thus we sequentially obtain
that $y_1,\ldots,y_k\in Y'$. Since $|Y'|<|Y|$, we may apply the inductive
assumption.

If $Y'=Y$, then $Y\subset H_1$, hence there exists $n\in\N$ such that if
$\bar H$ is the group
generated by $Y,TY,\ldots,T^nY$, then $Y\subset T\bar H$. Then $\bar H$ is
a finitely generated group, $T^{-1}$ an aperiodic endomorphism of $\bar H$.
For this case Lemma is already proved.
\newline\mbox{\ }\hfill\rule{2mm}{2mm}

Applying Lemma to the set $Y=X-X$ we see that the mapping
$$
X^k\to X_{nk},\ \ (x_1,\ldots,x_k)\mapsto\sum^k_{l=1}T^{n(l-1)}x_l,
$$
is a bijection for all $k\in\N$ and for all $n\in\N$ sufficiently large.
This bijection induces an isomorphism of $\Mat(X^k)$ onto
$\Mat(X_{nk})$. Composing it with
$i_{X_{nk}}\colon\Mat(X_{nk})\to C^*(G,\omega)$
and identifying $\Mat(X^k)$ with $\Mat(X)^{\otimes k}$ we obtain a unital
completely positive mapping
$$
\sigma_{nk}\colon\Mat(X)^{\otimes k}\to C^*(G,\omega).
$$
Set
$b(n,k)_{i_1\ldots i_k}=\sigma_{nk}(a_{i_1}\otimes\ldots\otimes a_{i_k})$.
By definition \cite{CNT}, we obtain

$\displaystyle
{1\over k}H_\tau(\gamma,\alpha^n_T\circ\gamma,
                  \ldots,\alpha^{n(k-1)}_T\circ\gamma)\ge
$
$$
\ge{1\over k}\sum_{i_1,\ldots,i_k}\eta\tau(b(n,k)_{i_1\ldots i_k})+
 {1\over k}\sum^k_{l=1}\sum_{i_l}S\Bigl(\tau\Bigl(\gamma(\cdot)\Bigr),
  \tau\Bigl(\gamma(\cdot)\alpha^{-n(l-1)}_T(b(n,k)^{(l)}_{i_l})\Bigr)\Bigr),
$$
where $\displaystyle b(n,k)^{(l)}_{i_l}
 =\sum_{i_1,\ldots,\hat i_l,\ldots,i_k}b(n,k)_{i_1\ldots i_k}$.

If we denote by $\tau_Y$ the unique tracial state on $\Mat(Y)$, then
$\tau_Y=\tau\circ i_Y$, so that $\tau\circ\sigma_{nk}=\tau^{\otimes k}_X$,
whence
$$
\tau(b(n,k)_{i_1\ldots i_k})=\prod^k_{l=1}\tau_X(a_{i_l})=
 \prod^k_{l=1}\tau(b_{i_l}).
$$
So the first term in the inequality above is equal to $\sum_i\eta\tau(b_i)$,
and in order to prove Theorem it remains to show that
$$
||\alpha^{-n(l-1)}_T(b(n,k)^{(l)}_{i_l})-b_{i_l}||\mathop{\to}_{n\to\infty}0
$$
uniformly on $k,l\in\N$ ($l\le k$) and $i_l\in I$. Let $\theta_l$ be
the embedding of $\Mat(X)$ into $\Mat(X)^{\otimes k}$ defined by
$$
\theta_l(a)=\underbrace{1\otimes\ldots\otimes1}_{l-1}\otimes a\otimes
 \underbrace{1\otimes\ldots\otimes1}_{k-l}.
$$
Then $b(n,k)^{(l)}_{i_l}=(\sigma_{nk}\circ\theta_l)(a_{i_l})$. Thus we just
have to estimate
$$
||\alpha^{-n(l-1)}_T\circ\sigma_{nk}\circ\theta_l-i_X||.
$$
Using the facts that $\omega$ is bilinear and $T$-invariant we obtain

$\displaystyle
\sigma_{nk}(e_{x_1x_1}\otimes\ldots\otimes e_{x_{l-1}x_{l-1}}\otimes e_{xy}
 \otimes e_{x_{l+1}x_{l+1}}\otimes\ldots\otimes e_{x_kx_k})=
$
\begin{eqnarray*}
 &=&{1\over|X|^k}\bar\omega\left(T^{n(l-1)}(x-y),T^{n(l-1)}y+
    \sum^k_{i=1,i\ne l}T^{n(i-1)}x_i\right)u_{T^{n(l-1)}(x-y)}\\
 &=&\left(\prod^k_{i=1,i\ne l}
      {\bar\omega(x-y,T^{n(i-l)}x_i)\over|X|}\right)
     {\bar\omega(x-y,y)\over|X|}u_{T^{n(l-1)}(x-y)},\\
\end{eqnarray*}
so that

$\displaystyle
||(\alpha^{-n(l-1)}_T\circ\sigma_{nk}\circ\theta_l-i_X)(e_{xy})||=
$
\begin{eqnarray*}
 &=&{1\over|X|}
    \left|\sum_{x_1,\ldots,\hat x_l,\ldots,x_k}\left(\prod^k_{i=1,i\ne l}
      {\bar\omega(x-y,T^{n(i-l)}x_i)\over|X|}\right)-1\right|\\
 &=&{1\over|X|}
    \left|\prod^k_{i=1,i\ne l}\left({1\over|X|}\sum_{z\in X}
      \bar\omega(x-y,T^{n(i-l)}z)\right)-1\right|.
\end{eqnarray*}

We must show that the latter expression tends to zero as $n\to\infty$
uniformly on $k,l\in\N$ ($l\le k$). This follows from
$$
\sum_{n\in\Z}\left|1-{1\over|X|}\sum_{z\in X}\omega(x-y,T^nz)\right|<\infty.
$$
So the proof of Theorem 2 is complete.

\bigskip

\section{Classical case}
If $\omega\equiv1$, then $C^*(G,\omega)=C(\hat G)$, the algebra of continuous
functions on the dual group~$\hat G$. It is known that an automorphism $T$ of
$G$ is aperiodic iff the dual automorphism of $\hat G$ is ergodic. Thus we
obtain a classical Rohlin's result \cite{R} stating that ergodic
automorphisms of
compact abelian groups have completely positive entropy. Note that in this
case we have
$$
b(n,k)_{i_1\ldots i_k}=b_{i_1}\alpha^n_T(b_{i_2})\ldots
    \alpha^{n(k-1)}_T(b_{i_k}),
$$
so what is really necessary for the proof is Lemma above and the possibility
of approximating in mean measurable partitions of unit by partitions
consisting from trigonometric polynomials, which can be proved by elementary
methods without appealing to Voiculescu's completely positive mappings.

\bigskip

\begin{flushleft}
Institute for Low Temperature Physics \& Engineering\\
Lenin Ave 47\\
Kharkov 310164, Ukraine\\
neshveyev@ilt.kharkov.ua\\
\end{flushleft}

\end{document}